\documentclass[a4paper,12pt]{article}%
\usepackage{amsmath}
\usepackage{amsfonts}
\usepackage{amssymb}
\usepackage{graphicx}%
\newtheorem{theorem}{Theorem}

\newtheorem{lemma}[theorem]{Lemma}

\newtheorem{remark}[theorem]{Remark}

\begin{document}

\title{Convolution operators defined by singular measures on the motion group}
\author{Luca Brandolini
\and Giacomo Gigante
\and Sundaram Thangavelu
\and Giancarlo Travaglini}
\date{}
\maketitle

\begin{abstract}
This paper contains an $L^{p}$ improving result for convolution operators
defined by singular measures associated to hypersurfaces on the motion group.
This needs only mild geometric properties of the surfaces, and it extends
earlier results on Radon type transforms on $\mathbb{R}^{n}$. The proof relies
on the harmonic analysis on the motion group.

\end{abstract}

\section{Introduction}

The classical Radon transforms satisfy $L^{p}$ improving properties (see
\cite{OS}) and they are closely related to certain convolution operators
associated to singular measures (see e.g. \cite{Str}). The above results have
been extended in many ways, not necessarily related to convolution structures,
see e.g. \cite{CNSW}, \cite{S}, \cite{TW} and the references therein.

Our starting point is the following result, proved in \cite{RT}.

\begin{theorem}
\label{Thm1}Let $\Gamma$ be a convex compact curve in the plane and let
$\gamma$ be the arc-length measure of $\Gamma$. We identify $\theta\in\left[
0,2\pi\right]  $ with $e^{i\theta}\in S^{1}$ (the unit circle). Let
$\gamma_{\theta}$ be the rotated measure, i.e. $\int_{\mathbb{R}^{2}}f\left(
x\right)  d\gamma_{\theta}\left(  x\right)  =\int_{\mathbb{R}^{2}}f\left(
e^{i\theta}x\right)  d\gamma\left(  x\right)  $. Consider the operator $T$
defined by%
\[
Tf\left(  x,\theta\right)  =\left(  f\ast_{\mathbb{R}^{2}}\gamma_{\theta
}\right)  \left(  x\right)  \;,
\]
where \ $x\in\mathbb{R}^{2}$ and $\ast_{\mathbb{R}^{2}}$ denotes the
convolution in $\mathbb{R}^{2}$. Then%
\[
\left\Vert Tf\right\Vert _{L^{3}\left(  \mathbb{R}^{2}\times S^{1}\right)
}\leq c\left\Vert f\right\Vert _{L^{3/2}\left(  \mathbb{R}^{2}\right)  }\;.
\]

\end{theorem}

The proof of this theorem relies on an estimate for the average decay of the
Fourier transform $\widehat{\gamma}$ proved by A.N. Podkorytov in \cite{P}
(see also \cite{BIT}), which has been extended to several variables in
\cite{BHI}. The following statement is different from the one in \cite{BHI},
but it can be proved by a mild variation of the original argument.

\begin{theorem}
\label{BHI-Improved}Let $\Gamma$ be a compact convex submanifold of
codimension $1$ in $\mathbb{R}^{n}$ (i.e. $\Gamma$ can be seen as the graph of
a convex function defined in a convex domain in $\mathbb{R}^{n-1}$). Let
$\gamma=\chi\sigma$ where $\sigma$ is the surface measure on $\Gamma$ and
$\chi$ is a smooth cutoff supported in the interior of $\Gamma$. Then%
\[
\int_{S^{n-1}}\left\vert \widehat{\gamma}\left(  R\omega\right)  \right\vert
^{2}d\omega\leq cR^{-\left(  n-1\right)  },
\]
where $d\omega$ is the normalized measure on the unit sphere $S^{n-1}$.
Moreover the constant $c$ depends only on $\chi$ and the diameter of $\Gamma$.
\end{theorem}

The above theorem easily implies the following extension of Theorem \ref{Thm1}
(see \cite{BGT}). For $k\in SO(n)$ and $\gamma$ a measure on $\mathbb{R}^{n}$,
let $\gamma_{k}$ be defined by $\int_{\mathbb{R}^{n}}fd\gamma_{k}%
=\int_{\mathbb{R}^{n}}f(ky)d\gamma\left(  y\right)  $, so that $\widehat
{\gamma}_{k}(\xi)=\widehat{\gamma}(k^{-1}\xi)$.

\begin{theorem}
\label{BGTradon}Let $\Gamma$ be a compact convex submanifold of codimension
$1$ in $\mathbb{R}^{n}$ and let $\gamma=\chi\sigma$ where $\sigma$ is the
surface measure on $\Gamma$ and $\chi$ is a smooth cutoff function supported
in the interior of $\Gamma$. Consider the operator $T$ defined by%
\[
Tf\left(  x,k\right)  =\left(  f\ast_{\mathbb{R}^{n}}\gamma_{k}\right)
\left(  x\right)  \;,
\]
where $x\in\mathbb{R}^{n}$, $k\in SO(n)$ and $\ast_{\mathbb{R}^{n}}$denotes
the convolution in $\mathbb{R}^{n}$. Then%
\[
\left\Vert Tf\right\Vert _{L^{n+1}\left(  \mathbb{R}^{n}\times SO(n)\right)
}\leq c\left\Vert f\right\Vert _{L^{\left(  n+1\right)  /n}\left(
\mathbb{R}^{n}\right)  }\;.
\]

\end{theorem}

The operator $T\ $in Theorem \ref{BGTradon} can be seen as a convolution
operator on the motion group $M_{n}$, which is $\mathbb{R}^{n}\times SO(n)$
equipped with the group product $\left(  x,k\right)  \left(  y,h\right)
=\left(  x+ky,kh\right)  $ and unit $\left(  0,e\right)  $. Indeed the
convolution of two functions $F$ and $G$ on $M_{n}$ is defined by%
\[
\left(  F\ast_{M_{n}}G\right)  \left(  x,k\right)  =\int_{M_{n}}F\left(
x-kh^{-1}y,kh^{-1}\right)  G\left(  y,h\right)  dydh,
\]
where $dh$ is the Haar measure on $SO(n)$.

Note that if $F\left(  x,k\right)  =f\left(  x\right)  $ and $\mu$ denotes the
measure on $M_{n}$ defined by%
\[
\int_{M_{n}}G\left(  x,k\right)  d\mu\left(  x,k\right)  =\int_{M_{n}}G\left(
x,k\right)  d\gamma_{k}\left(  x\right)  dk
\]
we have%
\begin{align}
F\ast_{M_{n}}\mu\left(  x,k\right)   &  =\int_{M_{n}}F\left(  x-kh^{-1}%
y,kh^{-1}\right)  d\mu\left(  y,h\right)  \label{asterisco}\\
&  =\int_{M_{n}}f\left(  x-kh^{-1}y\right)  d\gamma_{h}\left(  y\right)
dh\nonumber\\
&  =\int_{M_{n}}f\left(  x-ky\right)  d\gamma\left(  y\right)  dh\nonumber\\
&  =\int_{\mathbb{R}^{n}}f\left(  x-ky\right)  d\gamma\left(  y\right)
=\left(  f\ast_{\mathbb{R}^{n}}\gamma_{k}\right)  \left(  x\right)  .\nonumber
\end{align}

The above family $\left\{  \gamma_{k}\right\}  $ of hypersurfaces in
$\mathbb{R}^{n}$ turns out to be a manifold in $\mathbb{R}^{n}\times SO(n)$.
Indeed for any $k_{0}\in SO(n)$ the coset $\left\{  \left(  x,k_{0}\right)
:x\in\mathbb{R}^{n}\right\}  $ contains the $n-1$ dimensional manifold
$\Gamma_{k_{0}}$, i.e. the manifold $\Gamma$ rotated by $k_{0}$. The union of
the manifolds $\Gamma_{k_{0}}$ is a hypersurface $X$ in $\mathbb{R}^{n}\times
SO(n)$.

When $n=2$, the $\Gamma$'s are convex curves and their union can be seen as a
$2$-dimensional surface in $\mathbb{R}^{2}\times S^{1}$; the picture shows
this surface in the particular case $\Gamma\left(  t\right)  =\left(
t,t^{2}+1\right)  $, together with the plane $\theta=\pi$:%

\noindent\begin{center}
\begin{tabular}{ccc}
\includegraphics[width=12cm]{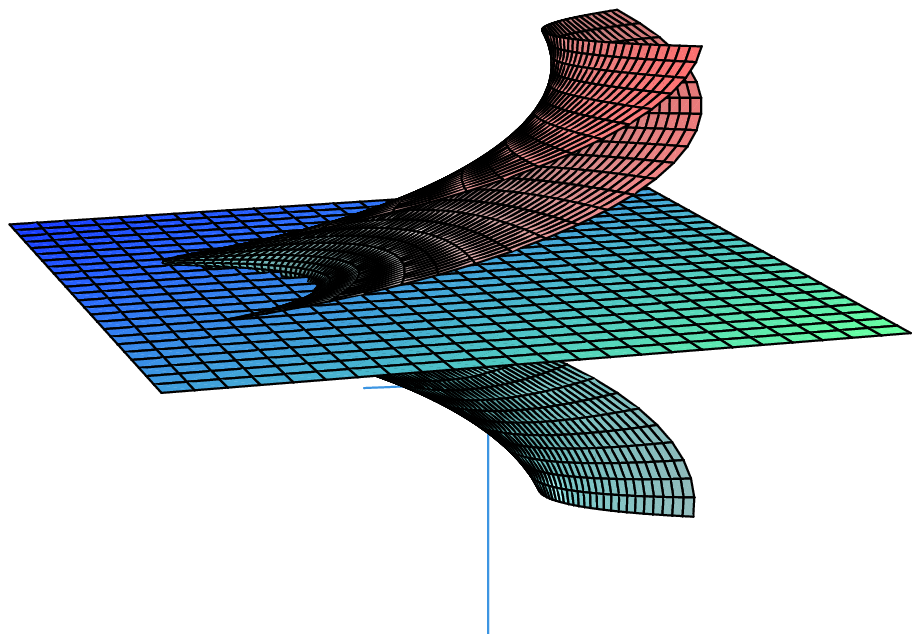} &\
\end{tabular}\\
\end{center}%

In this paper we want to replace the above manifold $X$ with a more general
manifold $Y$ in $\mathbb{R}^{n}\times SO(n)$, so that the action of $\Gamma$
as a convolution operator on $\mathbb{R}^{n}$ is averaged not only on
rotations, but on a wider family of transformations. In order to deal with
this more general setting it is natural to work in the Euclidean motion group
$M_{n}$ rather than in $\mathbb{R}^{n}\times SO(n)$ and take advantage of the
representation theory of $M_{n}$.

\section{Main result}

The following is our main result. By (\ref{asterisco}) it is an extension of
Theorem \ref{BGTradon}.

\begin{theorem}
\label{n}Let $n\geq2$ and let $Y$ be a $\mathcal{C}^{1}$ submanifold of
codimension $1$ in $M_{n}$. Assume that $Y$ can be locally represented as
$F\left(  x,k\right)  =0$ with $\nabla_{x}F\left(  x,k\right)  \neq0$. Assume
furthermore that for every $k_{0}\in SO(n)$ the intersection $Y\cap\left\{
\left(  x,k_{0}\right)  :x\in\mathbb{R}^{n}\right\}  $ is a convex
hypersurface\footnote{This means that the above intersection is the graph of a
convex function on a convex set (after choosing suitable coordinates in
$\mathbb{R}^{n-1}$).} in $\mathbb{R}^{n}$. Choose $\chi\in\mathcal{C}_{c}%
^{1}\left(  M_{n}\right)  $ and let $\mu$ be the measure on $M_{n}$ given by
$\int_{M_{n}}fd\mu=\int_{Y}f\chi d\sigma,$ where $\sigma$ is the surface
measure on $Y$. Then, if $Tf\left(  x,k\right)  =\left(  f\ast_{M_{n}}%
d\mu\right)  \left(  x,k\right)  $, we have%
\begin{equation}
\left\Vert Tf\right\Vert _{L^{n+1}\left(  M_{n}\right)  }\leq c_{n}\left\Vert
f\right\Vert _{L^{\left(  n+1\right)  /n}\left(  M_{n}\right)  } \label{bound}%
\end{equation}

\end{theorem}

\textit{Proof.} Without loss of generality, we may assume that $Y$ is the
graph of the function%
\[
x_{1}=\Phi\left(  x^{\prime},k\right)  ,
\]
where we use the notation $x^{\prime}=\left(  x_{2},\ldots,x_{n}\right)  $.
Thus
\[
\int_{M_{n}}fd\mu=\int_{\mathbb{R}^{n-1}}\int_{SO(n)}f\left(  \Phi\left(
x^{\prime},k\right)  ,x^{\prime},k\right)  \nu\left(  x^{\prime},k\right)
dkdx^{\prime},
\]
where $\nu$ is the product of $\chi$ by a Jacobian term. For every
$z\in\mathbb{C}$, let $i_{z}$ be the distribution on $\mathbb{R}$ defined by%
\[
\left\langle i_{z},\eta\right\rangle =\frac{1}{\Gamma\left(  z\right)  }%
\int_{0}^{+\infty}\eta\left(  t\right)  t^{z-1}dt\;.
\]
We define the family of distributions $\mu^{z}$ by%
\[
\mu^{z}=\mu\ast_{M_{n}}I_{z},
\]
where $I_{z}$ is the distribution defined by
\[
I_{z}\left(  x,k\right)  =i_{z}\left(  x_{1}\right)  \otimes\delta_{0}\left(
x_{2}\right)  \otimes\cdots\otimes\delta_{0}\left(  x_{n}\right)
\otimes\delta_{e}\left(  k\right)  .
\]
For any $k\in SO(n)$ define the measure $\mu_{k}$ on $\mathbb{R}^{n}$ by%
\[
\int_{\mathbb{R}^{n}}gd\mu_{k}=\int_{\mathbb{R}^{n-1}}g\left(  \Phi\left(
x^{\prime},k\right)  ,x^{\prime}\right)  \nu\left(  x^{\prime},k\right)
dx^{\prime}.
\]
Then define the distribution $E_{z}$ on $\mathbb{R}^{n}$ by%
\[
E_{z}\left(  x\right)  =i_{z}\left(  x_{1}\right)  \otimes\delta_{0}\left(
x_{2}\right)  \otimes\cdots\otimes\delta_{0}\left(  x_{n}\right)
\]
and let $\mu_{k}^{z}=\mu_{k}\ast_{\mathbb{R}^{n}}E_{z}$. Then it can be easily
shown that%
\begin{align*}
\int_{M_{n}}f\left(  x,k\right)  d\mu\left(  x,k\right)   &  =\int_{SO\left(
n\right)  }\int_{\mathbb{R}^{n}}f\left(  x,k\right)  d\mu_{k}\left(  x\right)
dk\\
\left\langle \mu^{z},f\right\rangle _{M_{n}}  &  =\int_{SO\left(  n\right)
}\left\langle \mu_{k}^{z},f\left(  \cdot,k\right)  \right\rangle
_{\mathbb{R}^{n}}dk.
\end{align*}
We introduce the analytic family of operators%
\[
T^{z}f=f\ast\mu^{z}.
\]
Then the proof follows from Stein's complex interpolation theorem and the
following result.

\begin{lemma}
For every real $s$ we have%
\begin{equation}
T^{1+is}:L^{1}\left(  M_{n}\right)  \longrightarrow L^{\infty}\left(
M_{n}\right)  \;, \label{1-inf}%
\end{equation}%
\begin{equation}
T^{-\left(  n-1\right)  /2+is}:L^{2}\left(  M_{n}\right)  \longrightarrow
L^{2}\left(  M_{n}\right)  \;. \label{2-2}%
\end{equation}

\end{lemma}

\textit{Proof of the Lemma. }Let us prove (\ref{1-inf}) first. Indeed for
$g\in L^{1}\left(  M_{n}\right)  $%
\begin{align*}
&  \left\langle \mu^{1+is},g\right\rangle _{M_{n}}=\left\langle \mu\ast
_{M_{n}}I_{1+is},g\right\rangle _{M_{n}}=\left\langle I_{1+is},g\ast_{M_{n}%
}\widetilde{\mu}\right\rangle _{M_{n}}\\
&  =\frac{1}{\Gamma\left(  1+is\right)  }\int_{0}^{+\infty}\left(
g\ast_{M_{n}}\widetilde{\mu}\right)  \left(  x_{1},0,\ldots,0,e\right)
~x_{1}^{is}dx_{1}\\
&  =\frac{1}{\Gamma\left(  1+is\right)  }\int_{0}^{+\infty}x_{1}^{is}%
~\int_{M_{n}}g\left(  \left(  x_{1},0,\ldots,0,e\right)  \left(  y_{1}%
,\ldots,y_{n},k\right)  ^{-1}\right) \\
&  d\widetilde{\mu}\left(  y_{1},\dots,y_{n},k\right)  dx_{1}\\
&  =\frac{1}{\Gamma\left(  1+is\right)  }\int_{0}^{+\infty}x_{1}^{is}%
~\int_{M_{n}}g\left(  \left(  x_{1},0,\ldots,0,e\right)  \left(  y_{1}%
,\ldots,y_{n},k\right)  \right) \\
&  d\mu\left(  y_{1},\dots,y_{n},k\right)  dx_{1}\\
&  =\frac{1}{\Gamma\left(  1+is\right)  }\int_{0}^{+\infty}x_{1}^{is}%
~\int_{M_{n}}g\left(  x_{1}+y_{1},y_{2},\ldots,y_{n},k\right)  d\mu\left(
y_{1},\dots,y_{n},k\right)  dx_{1}\\
&  =\frac{1}{\Gamma\left(  1+is\right)  }\int_{0}^{+\infty}\int_{\mathbb{R}%
^{n-1}}\int_{SO(n)}g\left(  x_{1}+\Phi\left(  y^{\prime},k\right)  ,y^{\prime
},k\right)  x_{1}^{is}\nu\left(  y^{\prime},k\right)  dkdy^{\prime}dx_{1}%
\end{align*}
\textit{ }(where $\widetilde{\mu}$ is defined by $\int_{M_{n}}f\left(
y,k\right)  d\widetilde{\mu}\left(  y,k\right)  =\int_{M_{n}}f\left(  \left(
y,k\right)  ^{-1}\right)  d\mu\left(  y,k\right)  \ $).

The substitution $y_{1}=x_{1}+\Phi\left(  y^{\prime},k\right)  $, along with
the boundedness of $\nu$, immediately gives
\[
\left\vert \left\langle \mu^{1+is},g\right\rangle _{M_{n}}\right\vert \leq
c\left\Vert g\right\Vert _{L^{1}\left(  M_{n}\right)  },
\]
so that $\mu^{1+is}\in L^{\infty}\left(  M_{n}\right)  $. This proves
(\ref{1-inf}).

Now we turn to the proof of (\ref{2-2}). We need first to recall a few facts
from the representation theory of $M_{n}$.

The unitary dual $\widehat{M}_{n}$ ($n\geq2$) can be described in the
following way (here \cite{ST} is a reference for the representation theory of
$M_{n}$, see also \cite{SU}). Let $L=SO\left(  n-1\right)  $, considered as a
subgroup of $SO(n)$. For each $\sigma\in\widehat{L}$ realised on a Hilbert
space $V_{\sigma}$ of dimension $d_{\sigma}$ consider the space $L^{2}\left(
SO(n),\sigma\right)  $ consisting of functions $\varphi$ on $SO(n)$ taking
values in $\mathbb{C}^{d_{\sigma}\times d_{\sigma}}$, the space of $d_{\sigma
}\times d_{\sigma}$ complex matrices, satisfying the condition%
\[
\varphi\left(  \ell k\right)  =\sigma\left(  \ell\right)  \varphi\left(
k\right)  \;,\;\;\;\;\;\;\;\;\;\;\;\;\;\ell\in L\;,\;\;\;\;\;k\in SO(n)
\]
which are also square integrable on $SO(n)$:%
\[
\int_{SO(n)}\left\Vert \varphi\right\Vert ^{2}dk=\int_{SO(n)}%
\operatorname*{tr}\left(  \varphi\left(  k\right)  ^{\ast}\varphi\left(
k\right)  \right)  dk<\infty\;.
\]
Note that $L^{2}\left(  SO(n),\sigma\right)  $ is a Hilbert space under the
inner product%
\[
\left(  \varphi,\psi\right)  =\int_{SO(n)}\operatorname*{tr}\left(
\varphi\left(  k\right)  ^{\ast}\psi\left(  k\right)  \right)  dk\;.
\]
For each $\lambda>0$ and $\sigma\in\widehat{L}$ we define a representation
$\pi_{\lambda,\sigma}$ of $M_{n}$ on $L^{2}\left(  SO(n),\sigma\right)  $ as
follows. For $\varphi\in L^{2}\left(  SO(n),\sigma\right)  $ and $\left(
x,k\right)  \in M_{n}$ let%
\[
\pi_{\lambda,\sigma}\left(  x,k\right)  \varphi\left(  \ell\right)
=\exp\left(  2\pi i\lambda\ell^{-1}e_{1}\cdot x\right)  \varphi\left(  \ell
k\right)  \;,
\]
where $e_{1}=\left(  1,0,\ldots,0\right)  $ and $\ell\in SO(n)$. If
$\varphi_{j}\left(  k\right)  $ are the column vectors of $\varphi\in
L^{2}\left(  SO(n),\sigma\right)  $ then $\varphi_{j}\left(  \ell k\right)
=\sigma\left(  \ell\right)  \varphi_{j}\left(  k\right)  $ for all $\ell\in
L$. Therefore $L^{2}\left(  SO(n),\sigma\right)  $ can be written as a direct
sum of $d_{\sigma}$ copies of $H\left(  SO(n),\sigma\right)  $ which is
defined to be the space of square integrable $\varphi:SO(n)\rightarrow
\mathbb{C}^{d_{\sigma}}$ satisfying%
\[
\varphi\left(  \ell k\right)  =\sigma\left(  \ell\right)  \varphi\left(
k\right)  ,\;\;\;\;\;\;\;\;\;\;\;\;\;\ell\in L\;.
\]
It can be shown that $\pi_{\lambda,\sigma}$ restricted to $H\left(
SO(n),\sigma\right)  $ is an irreducible representation of $M_{n}$. Moreover,
any infinite dimensional irreducible unitary representation of $M_{n}$ is
unitarily equivalent to one and only one $\pi_{\lambda,\sigma}$. Finite
dimensional irreducible unitary representations of $SO(n)$ also yield
irreducible unitary representations of $M_{n}$. As they do not appear in the
Plancherel formula we neglect them. We remark that when $n=2$ the unitary dual
$\widehat{L}$ contains only the trivial representation.

Given $f\in L^{1}\left(  M_{n}\right)  \cap L^{2}\left(  M_{n}\right)  $ we
define the group Fourier transform of $f$ by%
\[
\pi_{\lambda,\sigma}\left(  f\right)  =\int_{M_{n}}f\left(  x,k\right)
\pi_{\lambda,\sigma}\left(  \left(  x,k\right)  ^{-1}\right)  dxdk\;.
\]
It can be shown that $\pi_{\lambda,\sigma}\left(  f\right)  $ is a
Hilbert-Schmidt operator on $H\left(  SO(n),\sigma\right)  $ and we have the
Plancherel formula%
\[
\sum_{\sigma\in\widehat{L}}\int_{0}^{+\infty}\left\Vert \pi_{\lambda,\sigma
}\left(  f\right)  \right\Vert _{HS}^{2}\lambda^{n-1}d\lambda=\omega_{n}%
\int_{M_{n}}\left\vert f\left(  x,k\right)  \right\vert ^{2}dxdk\;,
\]
where $\left\Vert \cdot\right\Vert _{HS}$ denotes the Hilbert-Schmidt norm.

Applying Plancherel formula to $T^{-\left(  n-1\right)  /2+is}f$ we get%
\begin{align*}
&  \left\Vert T^{-\left(  n-1\right)  /2+is}f\right\Vert _{L^{2}\left(
M_{n}\right)  }^{2}\\
&  =\left\Vert f\ast_{M_{n}}\mu^{-\left(  n-1\right)  /2+is}\right\Vert
_{L^{2}\left(  M_{n}\right)  }^{2}\\
&  =\omega_{n}\sum_{\sigma\in\widehat{L}}\int_{0}^{+\infty}\left\Vert
\pi_{\lambda,\sigma}\left(  \mu^{-\left(  n-1\right)  /2+is}\right)
\pi_{\lambda,\sigma}\left(  f\right)  \right\Vert _{HS}^{2}\lambda
^{n-1}d\lambda\\
&  \leq\omega_{n}\sum_{\sigma\in\widehat{L}}\int_{0}^{+\infty}\left\Vert
\pi_{\lambda,\sigma}\left(  \mu^{-\left(  n-1\right)  /2+is}\right)
\right\Vert _{OP}^{2}~\left\Vert \pi_{\lambda,\sigma}\left(  f\right)
\right\Vert _{HS}^{2}\lambda^{n-1}d\lambda\;,
\end{align*}
where $\left\Vert \cdot\right\Vert _{OP}$ is the operator norm on $H\left(
SO(n),\sigma\right)  $. We shall show below that
\begin{equation}
\left\Vert \pi_{\lambda,\sigma}\left(  \mu^{-\left(  n-1\right)
/2+is}\right)  \right\Vert _{OP}\leq c_{n} \label{toprove}%
\end{equation}
uniformly in $\lambda$ and $\sigma$, so that
\[
\left\Vert T^{-\left(  n-1\right)  /2+is}f\right\Vert _{L^{2}\left(
M_{n}\right)  }^{2}\leq c_{n}\omega_{n}\sum_{\sigma\in\widehat{L}}\int
_{0}^{+\infty}\left\Vert \pi_{\lambda,\sigma}\left(  f\right)  \right\Vert
_{HS}^{2}\lambda^{n-1}d\lambda=c_{n}\left\Vert f\right\Vert _{L^{2}\left(
M_{n}\right)  }^{2}.
\]

We now prove (\ref{toprove}). For $\varphi,\psi\in H\left(  SO(n),\sigma
\right)  $ we have%
\[
\pi_{\lambda,\sigma}\left(  \left(  x,k\right)  ^{-1}\right)  \varphi\left(
u\right)  =\exp\left(  -2\pi i\lambda u^{-1}e_{1}\cdot k^{-1}x\right)
\varphi\left(  uk^{-1}\right)  .
\]
Assume for a moment $\operatorname{Re}z>0$, then $\mu^{z}$ is a measure and%
\begin{align*}
&  \pi_{\lambda,\sigma}\left(  \mu^{z}\right)  \varphi\left(  u\right) \\
&  =\int_{M_{n}}\pi_{\lambda,\sigma}\left(  \left(  x,k\right)  ^{-1}\right)
\varphi\left(  u\right)  d\mu^{z}\left(  x,k\right) \\
&  =\int_{M_{n}}\exp\left(  -2\pi i\lambda u^{-1}e_{1}\cdot k^{-1}x\right)
\varphi\left(  uk^{-1}\right)  d\mu^{z}\left(  x,k\right)  .
\end{align*}
Therefore%
\begin{align*}
&  \left\langle \pi_{\lambda,\sigma}\left(  \mu^{z}\right)  \varphi
,\psi\right\rangle _{H\left(  SO(n),\sigma\right)  }\\
&  =\int_{SO(n)}\left\langle \pi_{\lambda,\sigma}\left(  \mu^{z}\right)
\varphi\left(  u\right)  ,\psi\left(  u\right)  \right\rangle _{\mathbb{C}%
^{d_{\sigma}}}~du\\
&  =\int_{SO(n)}\int_{M_{n}}\exp\left(  -2\pi i\lambda u^{-1}e_{1}\cdot
k^{-1}x\right)  \left\langle \varphi\left(  uk^{-1}\right)  ,\psi\left(
u\right)  \right\rangle _{\mathbb{C}^{d_{\sigma}}}d\mu^{z}\left(  x,k\right)
du\\
&  =\int_{SO(n)}\int_{SO(n)}\int_{\mathbb{R}^{n}}\exp\left(  -2\pi
i\lambda\,k\,u^{-1}e_{1}\cdot x\right)  \left\langle \varphi\left(
uk^{-1}\right)  ,\psi\left(  u\right)  \right\rangle _{\mathbb{C}^{d_{\sigma}%
}}\\
&  \times d\mu_{k}^{z}\left(  x\right)  dudk\\
&  =\int_{SO(n)}\int_{SO(n)}\widehat{\mu_{k}^{z}}\left(  \lambda
k\,u^{-1}e_{1}\right)  \left\langle \varphi\left(  uk^{-1}\right)
,\psi\left(  u\right)  \right\rangle _{\mathbb{C}^{d_{\sigma}}}dudk\\
&  =\int_{SO(n)}\int_{SO(n)}\widehat{\mu_{k}}\left(  \lambda k\,u^{-1}%
e_{1}\right)  ~\widehat{E_{z}}\left(  \lambda k\,u^{-1}e_{1}\right)
\left\langle \varphi\left(  uk^{-1}\right)  ,\psi\left(  u\right)
\right\rangle _{\mathbb{C}^{d_{\sigma}}}dudk.
\end{align*}
By analytic continuation, the equality
\begin{align*}
&  \left\langle \pi_{\lambda,\sigma}\left(  \mu^{z}\right)  \varphi
,\psi\right\rangle _{H\left(  SO(n),\sigma\right)  }\\
&  =\int_{SO(n)}\int_{SO(n)}\widehat{\mu_{k}}\left(  \lambda k\,u^{-1}%
e_{1}\right)  ~\widehat{E_{z}}\left(  \lambda k\,u^{-1}e_{1}\right)
\left\langle \varphi\left(  uk^{-1}\right)  ,\psi\left(  u\right)
\right\rangle _{\mathbb{C}^{d_{\sigma}}}dudk
\end{align*}
holds also for $z=-\left(  n-1\right)  /2+is$. By Cauchy-Schwarz inequality
\begin{align*}
&  \left\vert \int_{SO(n)}\widehat{\mu_{k}}\left(  \lambda k\,u^{-1}%
e_{1}\right)  ~\widehat{E}_{-\left(  n-1\right)  /2+is}\left(  \lambda
k\,u^{-1}e_{1}\right)  \left\langle \varphi\left(  uk^{-1}\right)
,\psi\left(  u\right)  \right\rangle _{\mathbb{C}^{d_{\sigma}}}du\right\vert
\\
&  \leq\left\Vert \left\langle \varphi\left(  uk^{-1}\right)  ,\psi\left(
u\right)  \right\rangle _{\mathbb{C}^{d_{\sigma}}}\right\Vert _{L^{2}\left(
SO(n),du\right)  }\\
&  \times\left\Vert \widehat{\mu_{k}}\left(  \lambda k\,u^{-1}e_{1}\right)
~\widehat{E}_{-\left(  n-1\right)  /2+is}\left(  \lambda k\,u^{-1}%
e_{1}\right)  \right\Vert _{L^{2}\left(  SO(n),du\right)  }%
\end{align*}

By \cite[Ch. 2]{GS} we know that%
\[
\left\vert \widehat{E}_{-\left(  n-1\right)  /2+is}\left(  \lambda
k\,u^{-1}e_{1}\right)  \right\vert \leq C\lambda^{\left(  n-1\right)  /2}.
\]
Now%
\begin{align*}
\widehat{\mu_{k}}\left(  \xi\right)   &  =\int_{\mathbb{R}^{n}}\exp\left(
-2\pi i\xi\cdot x\right)  d\mu_{k}\left(  x\right)  \\
&  =\int_{\mathbb{R}^{n-1}}\exp\left(  -2\pi i\xi\cdot\left(  \Phi\left(
x^{\prime},k\right)  ,x^{\prime}\right)  \right)  \nu\left(  x^{\prime
},k\right)  dx^{\prime}\\
&  =\int_{\mathbb{R}^{n-1}}\exp\left(  -2\pi i\xi\cdot\left(  \Phi\left(
x^{\prime},k\right)  ,x^{\prime}\right)  \right)  \frac{\nu\left(  x^{\prime
},k\right)  }{\sqrt{1+\left\vert \nabla_{x^{\prime}}\Phi\left(  x^{\prime
},k\right)  \right\vert ^{2}}}\\
&  \times\sqrt{1+\left\vert \nabla_{x^{\prime}}\Phi\left(  x^{\prime
},k\right)  \right\vert ^{2}}dx^{\prime}\\
&  =\int_{\mathbb{R}^{n}}\exp\left(  -2\pi i\xi\cdot x\right)  \frac
{\nu\left(  x^{\prime},k\right)  }{\sqrt{1+\left\vert \nabla_{x^{\prime}}%
\Phi\left(  x^{\prime},k\right)  \right\vert ^{2}}}d\zeta_{k}\left(  x\right)
,
\end{align*}
where $d\zeta_{k}$ is the surface measure of the convex hypersurface in
$\mathbb{R}^{n}$ given by the intersection $Y\cap\left\{  \left(  x,k\right)
:x\in\mathbb{R}^{n}\right\}  $. By Theorem \ref{BHI-Improved} we get%
\begin{align*}
&  \left\Vert \widehat{\mu_{k}}\left(  \lambda k\,u^{-1}e_{1}\right)
~\widehat{E}_{-\left(  n-1\right)  /2+is}\left(  \lambda k\,u^{-1}%
e_{1}\right)  \right\Vert _{L^{2}\left(  SO(n),du\right)  }^{2}\\
&  \leq c\lambda^{n-1}\int_{SO(n)}\left\vert \widehat{\mu_{k}}\left(  \lambda
ku^{-1}e_{1}\right)  \right\vert ^{2}du\leq c.
\end{align*}
To end the proof we observe that%
\begin{align*}
&  \int_{SO(n)}\left\Vert \left\langle \varphi\left(  uk\right)  ,\psi\left(
u\right)  \right\rangle _{\mathbb{C}^{d_{\sigma}}}\right\Vert _{L^{2}\left(
SO(n),du\right)  }dk\\
&  =\int_{SO(n)}\left\{  \int_{SO(n)}\left\vert \left\langle \varphi\left(
uk\right)  ,\psi\left(  u\right)  \right\rangle _{\mathbb{C}^{d_{\sigma}}%
}\right\vert ^{2}du\right\}  ^{1/2}dk\\
&  \leq\left\{  \int_{SO(n)}\int_{SO(n)}\left\vert \left\langle \varphi\left(
uk\right)  ,\psi\left(  u\right)  \right\rangle _{\mathbb{C}^{d_{\sigma}}%
}\right\vert ^{2}dudk\right\}  ^{1/2}\\
&  \leq\left\{  \int_{SO(n)}\int_{SO(n)}\left\vert \varphi\left(  uk\right)
\right\vert ^{2}\left\vert \psi\left(  u\right)  \right\vert ^{2}dudk\right\}
^{1/2}\;.
\end{align*}
By Fubini's theorem and the invariance of the Haar measure on $SO(n)$ we get%
\[
\int_{SO(n)}\left\Vert \left\langle \varphi\left(  uk\right)  ,\psi\left(
u\right)  \right\rangle _{\mathbb{C}^{d_{\sigma}}}\right\Vert _{L^{2}\left(
SO(n),du\right)  }dk\leq\left\Vert \varphi\right\Vert _{H\left(
SO(n),\sigma\right)  }\left\Vert \psi\right\Vert _{H\left(  SO(n),\sigma
\right)  }.
\]
This ends the proof of the Lemma. Hence Theorem \ref{n} is proved.

\begin{remark}
For functions on $M_{n}$ which are independent of the rotational variable,
i.e. for functions $F$ such that $F\left(  x,k\right)  =f\left(  x\right)  $,
Theorem \ref{n} can be obtained from Theorem \ref{BGTradon}. Indeed%
\begin{align*}
&  \left\{  \int_{SO\left(  n\right)  }\int_{\mathbb{R}^{n}}\left\vert
F\ast_{M_{n}}\mu\left(  x,k\right)  \right\vert ^{n+1}dxdk\right\}
^{1/\left(  n+1\right)  }\\
&  =\left\{  \int_{SO\left(  n\right)  }\int_{\mathbb{R}^{n}}\left\vert
\int_{SO\left(  n\right)  }\int_{\mathbb{R}^{n}}f\left(  x-k\tau^{-1}y\right)
d\mu_{\tau}\left(  y\right)  d\tau\right\vert ^{n+1}dxdk\right\}  ^{1/\left(
n+1\right)  }\\
&  \leq\int_{SO\left(  n\right)  }\left\{  \int_{SO\left(  n\right)  }%
\int_{\mathbb{R}^{n}}\left\vert \int_{\mathbb{R}^{n}}f\left(  x-k\tau
^{-1}y\right)  d\mu_{\tau}\left(  y\right)  \right\vert ^{n+1}dxdk\right\}
^{1/\left(  n+1\right)  }d\tau\\
&  =\int_{SO\left(  n\right)  }\left\{  \int_{SO\left(  n\right)  }%
\int_{\mathbb{R}^{n}}\left\vert \int_{\mathbb{R}^{n}}f\left(  x-\sigma
y\right)  d\mu_{\tau}\left(  y\right)  \right\vert ^{n+1}dxd\sigma\right\}
^{1/\left(  n+1\right)  }d\tau\\
&  =\int_{SO\left(  n\right)  }\left\{  \int_{SO\left(  n\right)  }%
\int_{\mathbb{R}^{n}}\left\vert f\ast_{\mathbb{R}^{n}}\mu_{\tau,\sigma}\left(
x\right)  \right\vert ^{n+1}dxd\sigma\right\}  ^{1/\left(  n+1\right)  }%
d\tau\\
&  \leq\int_{SO\left(  n\right)  }c\left\Vert f\right\Vert _{L^{\frac{n+1}{n}%
}\left(  \mathbb{R}^{n}\right)  }d\tau=c\left\Vert F\right\Vert _{L^{\frac
{n+1}{n}}\left(  M_{n}\right)  },
\end{align*}
where $\mu_{\tau,\sigma}$ denotes the measure $\mu_{\tau}$ rotated by $\sigma
$. This yields the following weaker version of Theorem \ref{n}. For a general
$F$ let $\widetilde{F}\left(  x,k\right)  =\sup_{\tau\in SO\left(  n\right)
}\left\vert F\left(  x,k\tau\right)  \right\vert $ then%
\[
\left\Vert F\ast\mu\right\Vert _{L^{n+1}\left(  M_{n}\right)  }\leq\left\Vert
\widetilde{F}\ast\mu\right\Vert _{L^{n+1}\left(  M_{n}\right)  }\leq
c\left\Vert \widetilde{F}\right\Vert _{L^{\frac{n+1}{n}}\left(  M_{n}\right)
}=c\left\Vert F\right\Vert _{L_{x}^{\frac{n+1}{n}}\left(  L_{\tau}^{\infty
}\left(  M_{n}\right)  \right)  }.
\]
The above seems to be the best we can get by using earlier results such as the
ones in \cite{BGT}.
\end{remark}

\begin{remark}
A familiar example (the characteristic function of a small ball) and the
previous remark can be used to show that the indices in (\ref{bound}) cannot
be improved.
\end{remark}

\begin{remark}
It is interesting to compare Theorem \ref{n} with Theorem 1.1 in \cite{RS3}
where it is shown that the $L^{p}$ improving property of a measure is related
to the fact that the supporting manifold generates the full group.
\end{remark}

\begin{remark}
The techniques in our paper are $L^{2}$ in nature and they seem to provide
only $L^{p}-L^{p^{\prime}}$ results. We do not know how to get mixed norm
estimates similar to the ones which have been proved in \cite{RT} through
certain $L^{r}$ estimates for the average decay of Fourier transforms (note
that in general these $L^{r}$ estimates cannot be obtained by interpolation
between $L^{2}$ and $L^{\infty}$, see e.g. \cite{BRT}).
\end{remark}

\bigskip

The authors wish to thank the referee for her/his suggestions and for pointing
out some references.

\bigskip

\noindent\textsc{Dipartimento di Ingegneria dell'Informazione e Metodi
Mate\-matici}

\noindent\textsc{Universit\`{a} di Bergamo, Viale Marconi 5}

\noindent\textsc{24044 Dalmine, Bergamo, Italy}

\noindent\texttt{luca.brandolini@unibg.it}

\bigskip

\noindent\textsc{Dipartimento di Ingegneria dell'Informazione e Metodi
Mate\-matici}

\noindent\textsc{Universit\`{a} di Bergamo, Viale Marconi 5}

\noindent\textsc{24044 Dalmine, Bergamo, Italy}

\noindent\texttt{giacomo.gigante@unibg.it}

\bigskip

\noindent\textsc{Department of Mathematics}

\noindent\textsc{Indian Institute of Science}

\noindent\textsc{Bangalore 560012, India}

\noindent\texttt{veluma@math.iisc.ernet.in}

\bigskip

\noindent\textsc{Dipartimento di Statistica, Edificio U7}

\noindent\textsc{Universit\`{a} di Milano-Bicocca}

\noindent\textsc{Via Bicocca degli Arcimboldi 8}

\noindent\textsc{20126 Milano, Italy}

\noindent\texttt{giancarlo.travaglini@unimib.it}

\end{document}